%% file: HRC_final.tex
\documentclass[10pt]{article}
\usepackage{times}
\usepackage{amsmath,amsxtra,amsthm,amssymb,makeidx,graphics,amsfonts}
\usepackage{graphicx}
\usepackage{latexsym,color}
\usepackage{ytableau}
\usepackage[utf8x]{inputenc}
\input{preamble}

\title{Hook restriction coefficients}
\author{    \normalsize{ \bf{Sridhar P. Narayanan}} \\
            \em{School  of Mathematics} \\
            \em{GITAM Institute of Technology} \\
            \em{Vishakapatnam}\\
            email: sridharp.narayanan@gmail.com  }

\begin{document}

\maketitle
\begin{abstract}
The permutation matrices form a subgroup of $\GL_n(\C)$ that is isomorphic to the symmetric group $S_n$. Let $r_{\mu\lambda}$ denote the multiplicity of the irreducible representation $V_\mu$ of $S_n$, corresponding to a partition $\mu$ of $n$, in the restriction of an irreducible polynomial  representation $W_\lambda(\C)$ of $\GL_n(\C)$, corresponding to a partition $\lambda$ with at most $n$ parts. Finding a combinatorial interpretation for $r_{\mu\lambda}$ remains an open problem in algebraic combinatorics, called the \emph{restriction problem}. 

We derive a new nonrecursive expression for a character polynomial called the \emph{Specht polynomial} and use it to find a combinatorial interpretation of $r_{\mu\lambda}$ when $\lambda$ is a hook-shaped partition.
\end{abstract}

\noi

\section{Introduction}
\label{sec:intro}

The irreducible polynomial representations of $\GL_n(\C)$ are indexed by the set of partitions with at most $n$ parts. We denote them by $\WW$, where $\lambda$ is a partition with at most $n$ parts. The irreducible representations of the \emph{symmetric group} $S_n$ are indexed by the set of partitions of size $n$. We denote them by $\VV$, where $\mu$ is a partition of size $n$.
The $n\times n$ permutation matrices form a subgroup of $\GLn$ that is isomorphic to the symmetric group $S_n$. When $\WW$ is restricted to this subgroup, it decomposes as
\begin{displaymath}
\res^{\GLn}_{S_n}\left( \WW \right)= \bigoplus_{\mu \text{ partition of } n}\VV^{\oplus r_{\mu\lambda}}.
\end{displaymath}
The multiplicities $r_{\mu\lambda}$ are called \emph{restriction coefficients}. Finding a combinatorial interpretation for these coefficients is a longstanding open problem in algebraic combinatorics, called the \emph{restriction problem}.

Progress on this problem has been sporadic, and explicit combinatorial interpretations are scant. Assaf and Speyer in \cite{assaf} and  Orellana and Zabrocki in \cite{ORELLANA2021107943} independently found a basis of the symmetric functions that is connected to these coefficients. Section 5 of \cite{ORELLANA2021107943} connects their basis (called the \emph{irreducible character basis}) to the polynomials in this paper. In \cite{ALCO_2022__5_5_1165_0} Orellana, Zabrocki, Salioli and Schilling used the representation theory of partition algebras towards this problem. Sriwongsa, Heaton and Willenbring  found conditions for the positivity of sums of these coefficients in \cite{ALCO_2021__4_2_189_0}. Mitchell Lee in \cite{lee2025} found an interpretation for $r_{\mu\lambda}$ when $\lambda$ has at most three columns, using the following identity discovered by Littlewood in \cite{MR95209}:
$$s_\lambda[1+s_{(1)}+s_{(2)}+\dotsb]= \sum_{\mu}r_{\mu \lambda}s_\mu,$$
where $s_\lambda$ is the \emph{Schur polynomial} corresponding to a partition $\lambda$ (see \cite{MR3443860}) and $s_\lambda[f]$ denotes the  \emph{plethysm} of the polynomial $f$ into $s_\lambda$ (see \cite{MR1676282}). Some other important papers around these coefficients are \cite{10.1063/1.1666459, harman, 10.1063/1.1666632, MR95209, Scharf_1997}.


We \footnote{Amritanshu Prasad, Digjoy Paul, Shraddha Srivastava} use \emph{character polynomials} towards this problem. These are universal formulas that encode the characters of certain sequences of representations of the symmetric groups. This elementary yet promising approach has given us positivity conditions for $r_{(n)\lambda}$ (for the \emph{trivial representation} $V_{(n)}$) when $\lambda$ has either two rows, two columns or is hook-shaped in \cite{ALCO_2021__4_4_703_0}, and explicit combinatorial interpretations for the \emph{sign representation} $V_{(1^n)}$ when $\lambda$ has two columns or is hook-shaped in \cite{narayanan2022restriction}.

This paper extends those results to all $r_{\mu\lambda}$ when $\lambda$ is hook-shaped, which we call the \emph{hook restriction coefficients}. We begin by defining and recalling the basic concepts required for this paper in Section \ref{sec:pre}. In Section \ref{sec:Specht} we develop a nonrecursive formula for the Specht polynomial from a recursive formula given by Garsia and Goupil in \cite{MR2576382}. This formula, stated in Theorem \ref{th:specht}, may be of general interest in other problems involving FI-modules (see \cite{MR3357185}). We use it to express the restriction coefficient for a hook-shaped partition as a signed sum of coefficients in a generating function. Two sign-reversing involutions in Section \ref{sec:invol} eliminate negative terms, giving a final interpretation as the cardinality of a subset of \emph{supertableaus} in Theorem \ref{th:hrc}.
	\section{Preliminaries}
	\label{sec:pre}

This section recalls basic facts about compositions, partitions and tableau, and summarises the representation theory of symmetric groups and character polynomials. The reader is referred to \cite{fulton1991representation} for a comprehensive exposition of representation theory, to \cite{james_1984} specifically for symmetric groups, and to the work of Church, Ellenberg and Farb in \cite{MR3357185} on FI-modules and character polynomials.

A \defn{composition} is a finite sequence of \emph{nonnegative} integers. The \emph{size} of a composition is the sum of its parts, while the \emph{length} of a composition is the number of parts. We denote the size of a composition $\alpha$ by $|\alpha|$, and its length by $l(\lambda)$. We also use $\alpha \vDash n$ to denote that $\alpha$ is a composition of $n$. 

We depict compositions using \emph{Young diagrams}. The Young diagram of a composition $\alpha=(\alpha_1,\dotsb, \alpha_s)$ is a left-justified array of cells, with the $i$th row containing $\alpha_i$ cells, as in Figure \ref{fig:yngdcomp}. We say that $(i,j) \in \alpha$ if the Young diagram of $\alpha$ contains a cell in that position.
\begin{figure}[h]
\centering	
\begin{ytableau}
{\color{white}1} & {\color{white}1}&{\color{white}1} \\
{\color{white}1} &{\color{white}1} &{\color{white}1} &{\color{white}1}&{\color{white}1} &{\color{white}1}\\
{\color{white}1} &{\color{white}1}\\
{\color{white}1} &{\color{white}1} &{\color{white}1}\\
\end{ytableau}
\caption{The Young diagram for the composition $(3,6,2,3)$.}
\label{fig:yngdcomp}
\end{figure}
A \defn{partition} is a finite, \emph{nonincreasing} sequence of  \emph{positive} integers. The size and length of partitions is defined as for compositions. Let $\Par$ denote the set of all partitions, and let $\lambda \vdash n$ denote that $\lambda$ is a partition of size $n$. We often use the \emph{exponential notation} to represent a partition $\lambda$ as the formal sequence $1^{a_1}2^{a_2}\dotsb$, where $a_i$ is the number of parts of $\lambda$ equal to $i$. This paper deals with \emph{hook-shaped} partitions, which are of the form $(a+1,\underbrace{1\cdots 1}_{\text{$b$ times}})$, for $a,b\geq 0$. We employ the \emph{Frobenius notation} to depict this partition as $(a|b)$. 

Let $\mathcal{A}$ be a totally-ordered set. A \defn{tableau} $T$ of \emph{shape} $\alpha \vDash n$ in the \emph{alphabet} $\mathcal{A}$ is a filling of the Young diagram of $\alpha$ with elements of $\mathcal{A}$. We write $(i,j)\in T$ if $(i,j)\in \alpha$, and let $T(i,j)$ signify the entry in the cell $(i,j)$. We will define the tableaus required for this paper in Section \ref{sec:invol}

Let $S_n$ denote the $n$th symmetric group. The complex irreducible representations of $S_n$, called \emph{Specht modules}, are indexed by partitions of $n$. We denote them by $V_\mu$ and their characters by $\chi_\mu$ for a partition $\mu \vdash n$. 
\begin{example}
The representation corresponding to $(n)\vdash n$ is called the \emph{trivial representation}. Its character is the constant $\ch_{(n)}=1$.

The representation corresponding to $(1^n) \vdash n$ is called the \defn{sign representation}. Its character is given by
\begin{align}
\label{eq:sign}
\sgn(\sigma)= \prod_{i \geq 1}(-1)^{(i-1)x_i},
\end{align}  
where $x_i$ is the number of $i$-cycles in $\sigma$. 
\end{example}

Let $\GLn$ denote the $n$th general linear group. A representation $\rho:\GLn \rightarrow \GL_m(\mathbb{C})$ is called \emph{polynomial} if every entry of $\rho(A)$ is a polynomial in the entries of $A$, for every $A \in \GLn$. 
\begin{example}
Let $V=\Sym^{2}(\C^2)$ be the \emph{second symmetric tensor space}. Let $\{e_1,e_2\}$ be a basis of $\C^2$. Then the set $\{e_1^2,e_1e_2,e_2^2\}$ is a basis of $\Sym^2(\C^n)$, where $e_1e_2=e_2e_1$. The representation $\rho:\GLn\rightarrow \text{End}\left(\Sym^2(\C^n)\right)$ is defined on a matrix $A=\begin{pmatrix}a&b\\c&d\end{pmatrix}$ as 
\begin{displaymath}
\rho(A)=\begin{pmatrix}a^2 & ab & b^2\\2ac & ad+bc & 2bd\\c^2 & cd & d^2 \end{pmatrix}.
\end{displaymath}
Every entry of $\rho(A)$ is a polynomial in the entries of $A$.
\end{example} 
The irreducible polynomial representations of $\GLn$ are indexed by the set of partitions with \emph{at most} $n$ parts. Let $\WW$ denote the irreducible polynomial representation of $\GLn$ corresponding to a partition $\lambda$ with at most $n$ parts. The restriction of $\WW$ to the subgroup of permutation matrices, denoted $\res^{\GL_n}_{S_n}\left( \WW \right)$, is a representation of $S_n$. It therefore is the direct sum of Specht modules
\begin{displaymath}
\res^{\GLn}_{S_n}\left( \WW \right)= \bigoplus_{\mu \text{ partition of } n}\VV^{\oplus r_{\mu\lambda}}.
\end{displaymath}
We aim to find an interpretation for $r_{\mu\lambda}$ when $\lambda$ is hook-shaped.

An inner product is defined on $\mathcal{C}(n):=\{f:S_n\rightarrow \mathbb{C}|\text{ $f$ is \emph{constant} on each conjugacy class}\}$ as
\begin{align*}
\langle f,g \rangle_n = \frac{\sum_{\sigma \in S_n} f(\sigma) \overline{g(\sigma)}}{n!}.
\end{align*}
Since $f$ and $g$ are class-functions, the inner product above may be written as
\begin{align}
\label{eq:mom}
\langle f,g \rangle_n = \sum_{\alpha \vdash n} \frac{f(\alpha) \overline{g(\alpha)}}{z_{\alpha}},
\end{align}
where $f(\alpha), g(\alpha)$ are the values of $f$ and $g$ on the conjugacy class indexed by $\alpha=1^{a_1}2^{a_2}\dotsb\vdash n$, and $z_\alpha= \prod_{i\geq 1}i^{a_i}a_i!$ is the size of the centraliser of the conjugacy class indexed by $\alpha$. 
We abbreviate $\langle f, \ch_{(n)}\rangle_n$ to $\langle f\rangle_n$. 

The \emph{irreducible characters} $\{\ch_\mu|\mu\vdash n\}$ form an orthonormal basis of $\mathcal{C}(n)$. Every $f\in \mathcal{C}(n)$ can be written uniquely as $\sum_{a_\mu} a_\mu \ch_\mu$, where $a_\mu \in \mathbb{C}$. When $f$ is the character of a representation, the $a_\mu$ are nonnegative integers, that count the multiplicity of $V_\mu$ in the representation corresponding to $f$. In particular, when $f$ is the character of $\res^{\GL_n}_{S_n}\left( \WW \right)$ then $a_\mu=r_{\mu\lambda}$.

For each $i\geq 1$ and $n\geq 1$, let $x_i:S_n\rightarrow \mathbb{N}$ be defined on $\sigma \in S_n$ by 
\begin{displaymath}
x_i(\sigma)= \text{number of $i$-cycles of $\sigma$}.
\end{displaymath}

Polynomials in $\mathbb{C}[x_1,x_2,\dotsb]$ are called \emph{character polynomials}. We may define the evaluation of these polynomials at a permutation $\sigma$ by $p(\sigma)= p(x_1(\sigma),x_2(\sigma),\dotsb)$. A sequence $V_n$ of represenations of $S_n$ is said to be \emph{eventually polynomial} if there exist integers $n$ and $r$ and $p\in \mathbb{C}[x_1,x_2,\dotsb]$ such that
\begin{displaymath}
\chi_{V_n}(\sigma)= p(\sigma),
\end{displaymath}
for all $\sigma \in S_n$ for $n\geq N$. The polynomial $p$ is called the \textbf{character polynomial} of the sequence $(V_n)$. Character polynomials act as universal character formulae for sequences of representations. 

Given the partition $\mu=(\mu_1,\dotsb,\mu_d)$, define the padded partition $\mu[n]=(n-|\mu|,\mu_1,\dotsb,\mu_d)$ for all $n\geq |\mu|+\mu_1$. The sequence of irreducible representations $(V_{\mu[n]})$ is eventually polynomial (for $n\geq |\mu|+\mu_1$). We call the associated character polynomial a \textbf{Specht polynomial}, and denote it by $q_\mu$. Formulas for the Specht polynomials occur in \cite{MR2576382} and \cite{MR3443860}. We use the latter to find a derive a new formula for Specht polynomials in Section \ref{sec:Specht}.
In \cite{ALCO_2021__4_4_703_0} we find an expression for the character polynomial of $(\res^{\GL_n}_{S_n}\left( \WW \right)$, and in particular for the polynomials $H_k$ and $E_\ell$ corresponding to the sequences $(\text{Sym}^k(\C^n))$ and $(\bigwedge^l(\C^n))$ respectively. When $\lambda=(k|\ell)$, the Pieri rule (see \cite{MR3443860}) allows us to express
\begin{displaymath}
S_{(k|\ell)}= \sum_{j=0}^\ell(-1)^j H_{k+1+j}E_{\ell-j}.
\end{displaymath} 

We aim to find generating functions for $r_{\mu (k|\ell)}=\langle q_\mu,S_{(k|\ell)}\rangle_n$, and use combinatorial techniques thereafter to simplify to a positive sum.

			\section{The Specht polynomial}
			\label{sec:Specht}
Fix a partition $\mu$ and define $\mu[n]=(n-|\mu|,\mu_1,\dotsb,\mu_{d})$ for $n\geq |\mu|+\mu_1$. The sequence of Specht modules $(V_{\mu[n]})_{n\geq |\mu|+\mu_1}$ is eventually polynomial, with its character polynomial given by the following formula by Garsia and Goupil in \cite[Theorem 1.1]{MR2576382}:
\begin{align}
\label{eq:qrec}
q_{\mu}(x_1,x_2,\dotsb)= \sum_{\alpha \vdash |\mu|}\frac{q_{r(\mu)}(\alpha)}{z_\alpha}\prod_i (ix_i-1)^{a_i}\downarrow,
\end{align} 
where $\alpha=1^{a_1}2^{a_2}\dotsb$, and $r(\mu)=(\mu_2,\dotsb,\mu_{d})$ and the \emph{umbral operator} $\downarrow$ is defined as $x_i^{b}\downarrow=b!{x_i \choose b}$.
 
\begin{remark}[\defn{Notation}]
For partitions $\alpha=1^{a_1}2^{a_2}\dotsb$ and $\beta=1^{b_1}2^{b_2}\dotsb$, we let ${\alpha \choose \beta}:= \prod_{i\geq 1} {{a_i}\choose {b_i}}$ and ${\mathbf{x}\choose \beta}:= \prod_{i\geq 1} {{x_i}\choose {b_i}}$ for $\mathbf{x}=1^{x_1}2^{x_2}\dotsb$.
\end{remark}
\begin{prop}
\label{prop:Qrec}
For a partition $\beta$, let $B(\mu,\beta):= (-1)^{\sum x_i}{\mathbf x \choose \beta}q_{r(\mu)}(\mathbf x)\in \mathbb{C}[x_1,x_2,\dotsb]$.  Then 
\begin{displaymath}
q_\mu(x_1,x_2,\dotsb)=\sum_{\beta} (-1)^{l(\beta)}z_\beta{\mathbf{x} \choose \beta}\langle B(\mu,\beta)\rangle_{|\mu|},
\end{displaymath}

\end{prop}
\begin{proof}
From Equation \eqref{eq:qrec} we have
\begin{align*}
q_{\mu}(x_1,x_2,\dotsb)&= \sum_{\alpha \vdash |\mu|}\frac{q_{r(\mu)}(\alpha)}{z_\alpha}\prod_i (ix_i-1)^{a_i}\downarrow\\
&= \sum_{\alpha \vdash |\mu|} \frac{q_{r(\mu)}(\alpha)}{z_\alpha}\prod_{i\geq1} \sum_{b_i\leq a_i} (-1)^{a_i-b_i}{a_i \choose b_i}i^{b_i}x_i^{b_i}\downarrow\\
&=\sum_{\alpha \vdash |\mu|} \frac{q_{r(\mu)}(\alpha)}{z_\alpha}\prod_{i\geq 1} \sum_{b_i\leq a_i} (-1)^{a_i-b_i}{a_i \choose b_i}i^{b_i}b_i!{x_i\choose b_i}\\
&= \prod_{i \geq 1}\sum_{b_i\geq 0}(-1)^{b_i} i^{b_i}b_i!{x_i\choose b_i}\sum_{\alpha\vdash |\mu|}\frac{(-1)^{a_i}q_{r(\mu)}(\alpha)}{z_\alpha}{a_i \choose b_i}\\
&=\sum_{\beta \in \Par}(-1)^{l(\beta)}z_\beta {\mathbf x \choose \beta} \sum_{\alpha \vdash |\mu|} \frac{(-1)^{l(\alpha)}q_{r(\mu)}(\alpha)}{z_\alpha}{\alpha \choose \beta},
\end{align*}
where $\sum_{\alpha \vdash |\mu|} \frac{(-1)^{l(\alpha)}q_{r(\mu)}(\alpha)}{z_\alpha}{\alpha \choose \beta}= \langle B(\mu,\beta)\rangle_{|\mu|}$.
\end{proof}
We unwind the recursion in Equation \eqref{eq:qrec} by finding a generating function for $\langle B(\mu,\beta)\rangle_{|\mu|}$. We will make frequent use of the following identities to prove results in this section and the next:
\begin{gather*}
\tag A
\label{eq:exp}
e^x= \sum_{n\geq 0} \frac{x^n}{n!}\\
\tag B
\label{eq:log}
ln(1-x)= -\sum_{n\geq 1} \frac{x^n}{n}.
\end{gather*}

\begin{example}
\label{eg:Qmu1}
When $\mu=(\mu_1)$, then $B(\mu,\beta)=(-1)^{\sum x_i}{x \choose \beta}$, since $q_\emptyset=\ch_{(n)}=1$. In this case
\begin{align*}
\sum_{\mu \in \Par}\langle B(\mu,\beta) \rangle_{|\mu|}w^{|\mu|}&{=}\sum_{\mu \in \Par, \alpha \vdash |\mu|}\frac{(-1)^{l(\alpha)}w^{|\alpha|}}{z_{\alpha}}{\alpha \choose \beta}\\
&=\prod_{i\geq 1}\sum_{a_i\geq 0}\frac{(-1)^{a_i}w^{ia_i}}{i^{a_i}a_i!}{a_i \choose b_i}\\
&=\prod_{i\geq 1}\frac{(-1)^{b_i}w^{ib_i}}{b_i!i^{b_i}}\sum_{a_i\geq 0}\frac{(-1)^{a_i-b_i}w^{i(a_i-b_i)}}{i^{(a_i-b_i)}(a_i-b_i)!}\\ 
&= \frac{(-1)^{l(\beta)}w^{|\beta|}}{z_\beta}\prod_{i\geq 1}\sum_{c_i\geq 0} \frac{(-w^i)^{c_i}}{i^{c_i}c_i!}, \quad\quad \text{where $c_i=a_i-b_i$}\\
&\overset{\ref{eq:exp}}{=} \frac{(-1)^{l(\beta)}w^{|\beta|}}{z_\beta}\prod_{i \geq 1}\text{exp}\left(\frac{-w^i}{i}\right)\\
&\overset{\ref{eq:log}}{=} \frac{(-1)^{l(\beta)}w^{|\beta|}}{z_\beta}\left(1-w\right)\\
\end{align*}
\end{example}

\begin{definition}
\label{def:sum}
For a partition $\mu=(\mu_1,\dotsb,\mu_d)$, define $(\sum \mu)\in \Par$ by $$(\sum \mu)_i= \sum_{j=i}^d\mu_j.$$
\end{definition}
\begin{definition}
\label{def:Pd}
Define $P_d(w_1,\dotsb,w_d)$ by the recurrence
\begin{align*}
P_d(w_1,\dotsb,w_d)&= (1-w_1)(1-w_1w_2)\dotsb (1-w_1\dotsb w_d) P_{d-1}(w_2,\dotsc,w_d),\\
P_1(w_d)&=(1-w_d).
\end{align*} 
\end{definition}
\begin{prop}
\label{prop:Qform}
Fix a partition $\mu=(\mu_1,\dotsb,\mu_d)$ with at most $d$ parts, and a partition $\brac{1}$. We have:
\begin{displaymath}
\sum_{\mu: l(\mu)\leq d}\langle B(\mu,\brac{1})\rangle_{|\mu|}\mbf w^{\sum \mu}= \frac{(-1)^{l({\brac{1}})}}{z_{{\brac{1}}}}P_d(\mbf w)\sum_{\brac{2},\dotsc,\brac{d} \in \text{Par}}{{\brac{1}} \choose \brac{2}}\dotsb{\brac{d-1} \choose \brac{d}}\mbf w^{|\mathbf{\beta}|}, 
\end{displaymath}
where $\mbf{w}^{\sum \mu}=w_1^{|\mu|}w_2^{(\sum \mu)_2}\dotsb w_d^{\mu_d}$ and $\mbf{w}^{|\mathbf{\beta}|}=w_1^{|\brac{1}|}\dotsb w_d^{|\brac{d}|}$.
\end{prop}
\begin{proof}
We prove this by induction on $d$. Example \ref{eg:Qmu1} handles the basic step. Assume the result is true for all partitions with at most $d-1$ parts, in particular for $r(\mu)$. Using Proposition \ref{prop:Qrec}:
\begin{align*}
\sum_{\mu: l(\mu)\leq d} B(\mu,\brac{1})\overline{\mbf w}^{\sum r(\mu)}&=\sum_{\mu: l(\mu)\leq d} (-1)^{\sum x_i}{x \choose \brac{1}}q_{r(\mu)}(x)\overline{\mbf w}^{\sum r(\mu)}\\
&=(-1)^{\sum x_i}{x \choose \brac{1}}\sum_{\brac{2}} (-1)^{l(\brac{2})}z_{\brac{2}}{x \choose \brac{2}}\sum_{\mu:\mu=\nu \text{ or }r(\mu)=\nu} \sum_{\nu:l(\nu)\leq d-1}\langle B(\nu,\brac{2})\rangle_{|\nu|}\overline{\mbf w}^{\sum \nu},
\end{align*}
where $\overline{\mbf w}=(w_2,\dotsb, w_d)$. Using the inductive hypothesis on $\langle B(\nu,\brac{2})\rangle_{|\nu|}$, we have:
\begin{displaymath}
\sum_{\nu| l(\nu)\leq d-1}\langle B(\nu,\brac{2})\rangle_{|\nu|}\overline{\mbf w}^{\sum \nu}= \frac{(-1)^{l({\brac{2}})}}{z_{{\brac{2}}}}P_{d-1}(\overline{\mbf w})\sum_{\brac{3},\dotsc,\brac{d} \in \text{Par}}{{\brac{2}} \choose \brac{3}}\dotsb{\brac{d-1} \choose \brac{d}}\mbf {w}^{|\overline{\mathbf{\beta}}|}, 
\end{displaymath}

where $|\overline{\mbf \beta}|=(|\brac{2}|,\dotsb,|\brac{d}|)$. Thus
\begin{align*}
\sum_{\mu}\langle B(\mu,\brac{1})\rangle_{|\mu|}\mbf{w}^{\sum \mu}=P_{d-1}(\overline{\mbf{w}}) \sum_{\brac{2},\dotsc,\brac{d}}\sum_{\alpha \vdash |\mu|}w_1^{|\mu|} \frac{(-1)^{l(\alpha)}}{z_\alpha} {\alpha \choose \brac{1}}{\alpha \choose \brac{2}}{\brac{2} \choose \brac{3}}\dotsb {\brac{d-1} \choose \brac{d}}\overline{\mbf{w}}^{|\overline{\beta}|}.
\end{align*}
Let $\beta^{(j)}= 1^{b^{(j)}_1}2^{b^{(j)}_2}\dotsb$ in the exponential notation. Then $$\sum_{\mu}\langle B(\mu,\brac{1})\rangle_{|\mu|}\mbf{w}^{\sum \mu}=P_{d-1}(\overline{\mbf{w}})\prod_{i\geq 1}T_i,$$
where
\begin{align*}
T_i&= \sum_{a_i\geq 0}\sum_{\brac[b]{2}_i,\dotsb,\brac[b]{d}_i} \frac{(-1)^{a_i}}{i^{a_i}a_i!}{a_i \choose \brac[b]{1}_i}{a_i \choose \brac[b]{2}_i}{\brac[b]{2}_i \choose \brac[b]{3}_i}\dotsb {\brac[b]{d-1}_i \choose \brac[b]{d}_i} w_1^{ia_i}w_2^{i\brac[b]{2}_i}\dotsb w_d^{i\brac[b]{d}_i},
\end{align*}
where each $\brac[b]{j}_i\geq 0$. Separating out terms belonging to $\brac{1}$:
\begin{align*}
T_i&= \frac{(-1)^{\brac[b]{1}_i}w_1^{i\brac[b]{1}_i}}{i^{\brac[b]{1}_i}\brac[b]{1}_i!}\sum_{c_i\geq 0}{a_i \choose \brac[b]{1}_i} \sum_{\brac[b]{2}_i,\dotsb,\brac[b]{d}_i}\frac{(-1)^{c_i}}{i^{c_i}c_i!}{a_i \choose \brac[b]{2}_i}{\brac[b]{2}_i \choose \brac[b]{3}_i}\dotsb {\brac[b]{d-1}_i \choose \brac[b]{d}_i} w_1^{ic_i}w_2^{ib^{(2)}_i}\dotsb w_d^{ib^{(d)}_i},
\end{align*}
with $c_i=a_i-\brac[b]{1}_i$.

Using Equation \eqref{eq:exp}, we have $T_i= W_i^{\brac[b]{1}_i}\text{exp}\left(\frac{-W_i}{i}\right)$,
where $W_i= w_1^i(1+w_2^i(1+w_3^i(1+\dotsc +w_{d-1}^i(1+w_d^i))\dotsc)$. Thus

\begin{align*}
\sum_{\mu}\langle B(\mu,\brac{1})\rangle_{|\mu|}\mbf{w}^{|\sum \mu|}&=\frac{(-1)^{l(\brac{1})}}{z_{\brac{1}}}P_{d-1}(\overline{w})\prod_{i\geq 1}\text{exp}\left(\frac{-W_i}{i}\right)W_i^{\brac[b]{1}_i}.
\end{align*}
Using Equation \eqref{eq:log}:
\begin{align*}
\prod_{i\geq 1}\text{exp}\left(\frac{-W_i}{i}\right)&= (1-w_1)(1-w_1w_2)\dotsb (1-w_1\dotsb w_d),
\end{align*}
and through repeated applications of the binomial theorem:
\begin{align*}
\prod_{i\geq 1} W_i^{\brac[b]{1}_i}&= \sum_{\brac{2},\dotsb, \brac{d}}{\brac{1} \choose \brac{2}}\dotsb {\brac{d-1} \choose \brac{d}}w_1^{|\brac{1}|}\dotsb w_d^{|\brac{d}|}
\end{align*}

Substituting back yields the desired expression.
\end{proof}

Propositions \ref{prop:Qrec} and \ref{prop:Qform} give us a {new formula} for the Specht polynomial. 
\begin{theorem}[\textbf{Specht polynomial}]
\label{th:specht}
The generating functions for Specht polynomials is
\begin{align}
\label{eq:spchar}
\sum_{\mu: l(\mu)\leq d}q_\mu(\mathbf x) \mbf{w}^{|\sum \mu|}=P_d(\mbf w)\sum_{\brac{1},\dotsc,\brac{d}\in \text{Par}} {\mathbf x \choose \brac{1}}{\brac{1} \choose \brac{2}}\dotsb {\brac{d-1}\choose \brac{d}}\mbf w^{|\mbf \beta|}.
\end{align}

\end{theorem}

		\section{Generating function for hook restriction coefficients}
		\label{sec:genfunc}
Let $\text{Sym}^k(\C^n)$ be the $k$-th symmetric tensor space of $\C^n$, and let $\bigwedge^\ell(\C^n)$ be the $\ell$-th alternating tensor space of $\C^n$. The sequences $(\text{Sym}^k(\C^n))_{n\geq 1}$ and $(\bigwedge^\ell(\C^n))_{n\geq 1}$ are eventually polynomial. Let $H_k(\mathbf{x})$ and $E_\ell(\mathbf{x})$ denote the {character polynomials} for $\text{Sym}^k(\C^n)$ and $\bigwedge^\ell(\C^n)$ respectively. The following formulas for them were found in \cite{ALCO_2021__4_4_703_0}:
\begin{align}
\label{eq:hkel}
H_k(\mathbf{x})= \sum_{\gamma \vdash k} \mch{\mathbf{x}}{\gamma},\nonumber\\
E_\ell(\mathbf{x})= \sum_{\nu \vdash \ell} (-1)^{|\nu|-l(\nu)}{\mathbf{x} \choose {\nu}},
\end{align}
where $\mch{x_i} {g_i}={{x_i+g_i-1}\choose {g_i}}$ is the number of multisets of $g_i$ elements drawn from a set with $x_i$ elements.
\begin{prop}
\label{lem:a}
Let $\mbf{t}=(t_1,t_2,\dotsb,t_{d})$. Then
\begin{displaymath}
\sum_{\substack{k,\ell,\mu\\n \geq |\mu|+\mu_1}}\langle q_{\mu}, H_{k}E_{\ell}\rangle_{n}u^kv^\ell \mbf{t}^{\mu[n]}= \Upsilon\frac{\prod_{j=1}^{d}\prod_{r\geq 0}(1+u^rvt_j)}{\prod_{j=1}^{d}\prod_{r\geq 0}(1-u^rt_j)},
\end{displaymath} 
where the sum is over partitions $\mu$ with at most $d-1$ parts and $\Upsilon:= \prod_{i>j}(1-\frac{t_i}{t_j})$.
\end{prop} 
\begin{proof}
Using Equation \eqref{eq:mom}, Theorem \ref{th:specht} and \eqref{eq:hkel}, we have 
\begin{align*}
\sum_{k,\ell,\mu}\langle q_{\mu}, H_{k}E_{\ell}\rangle_{n}w_0^{n}u^kv^\ell \mbf{w}^{\mbf{\sum \mu}}&=\\
P(\mbf w)\sum_{\alpha \vdash n}\frac{w_0^{|\alpha|}}{z_\alpha}\sum_{k,\ell,\mu} u^kv^\ell\sum_{\gamma \vdash k} \mch{\mathbf{\alpha}}{\gamma}  \sum_{\nu \vdash \ell}  (-1)^{|\nu|-l(\nu)}{\mathbf{\alpha} \choose {\nu}}& \sum_{\brac{1},\dotsc,\brac{d-1}\in \text{Par}} \mbf w^{|\mbf \beta|}{\alpha \choose \brac{1}}{\brac{1} \choose \brac{2}}\dotsb {\brac{d-2}\choose \brac{d-1}},
\end{align*}
where $\mbf{w}=(w_1,\dotsb,w_{d-1})$.
We will follow the notation $\alpha=1^{a_1},\dotsb$ and $\brac{1}=1^{\brac[b]{1}_1},\dotsb$ and $\gamma= 1^{g_1}\dotsb$ and $\nu=1^{c_1}\dotsb$ and let $\mbf{b}_i=(\brac[b]{1}_i,\dotsb, \brac[b]{d-1}_i)$. The expression above is the product
\begin{displaymath}
P_{d-1}(\mbf w)\prod_{i \geq 1}B_i,
\end{displaymath}
where $$B_i=\sum_{a_i\geq 0}\frac{w_0^{ia_i}}{i^{a_i}a_i!}\sum_{g_i\geq 0, c_i\geq 0}  \mch{\mathbf{a_i}}{g_i}u^{ig_i}v^{ic_i} (-1)^{(i-1)c_i}{a_i \choose c_i} \sum_{\brac{1},\dotsc,\brac{d-1}} \mbf w^{i\mbf{b}_i}{a_i \choose \brac[b]{1}_i}\dotsb {\brac[b]{d-2}_i\choose \brac[b]{d-1}_i}.$$

We then use the identities
\begin{align*}
(1+x)^n= \sum_{0\leq k \leq n} {n \choose k} x^k,\\
(1-x)^{-n}= \sum_{k\geq 0} \mch{n}{k} x^k
\end{align*}
to express $\sum_{\substack{k,\ell,\mu\\n\geq |\mu|+\mu_1}}\langle q_{\mu}, H_{k}E_{\ell}\rangle_{n}u^kv^\ell w_0^{n}\mbf{w}^{\mbf{\sum \mu}}$ as
\begin{displaymath}
P_{d-1}(\mbf w)\prod_{i \geq 1}\text{exp}\left[\frac{(W_i)(1-u^i)^{-1}(1-(-v)^{i})}{i}\right],
\end{displaymath}
where $W_i= w_0^i(1+w_1^i(1+w_2^i(1+w_3^i(1+\dotsc +w_{d-1}^i(1+w_{d-1}^i))\dotsc)$. 

Then
\begin{align*}
P_{d-1}(\mbf w)\prod_{i \geq 1}\text{exp}\left[\frac{W_i(1-u^i)^{-1}(1-(-v)^{i})}{i}\right]&=P_{d-1}(\mbf w)\prod_{i \geq 1}\text{exp}\left[\frac{W_i(\sum_{r\geq 0}u^{ir})(1-(-v)^{i})}{i}\right]\\
=P_{d-1}(\mbf w)\prod_{r\geq 0}&\text{exp}\left[\frac{\sum_{i \geq 1}W_iu^{ir}}{i}\right]\text{exp}\left[-\frac{\sum_{i \geq 1}W_i(-u^{r}v)^{i}}{i}\right]\\
&=P_{d-1}(\mbf w)\prod_{r,0\leq j\leq d-1}\frac{1+\brac[W]{j}u^rv}{1-\brac[W]{j}u^r},
\end{align*}
where $\brac[W]{j}= w_0w_1\dotsb w_j$.

In summation, we have proved:
\begin{displaymath}
\sum_{\substack{k,\ell,\mu\\n \geq |\mu|+\mu_1}}\langle q_{\mu}, H_{k}E_{\ell}\rangle_{n}w_0^{n}u^kv^\ell \mbf{w}^{\mbf{\sum \mu}}= P_{d-1}(\mbf w)\prod_{r\geq 0,0\leq j\leq d-1}\frac{1+\brac[W]{j}u^rv}{1-\brac[W]{j}u^r}.
\end{displaymath}
Substitute $w_0=t_1$ and $w_i= \frac{t_{i+1}}{t_i}$ for $1\leq i \leq d$, whereby each $\brac[W]{j}=t_j$ and $P_{d-1}(\mbf w)=\Upsilon(t_1,\dotsb,t_d)$, giving the desired expression. 
\end{proof}
We know by the Pierei rule that 
\begin{align*}
S_{(k|\ell)}(\mathbf{x})= \sum_{j=0}^{\ell} (-1)^{j}H_{k+1+j}E_{l-j}.
\end{align*}

We can express $r_{\mu[n](k|\ell)}=\langle q_\mu,S_{(k|\ell)}\rangle_n$ as an alternating sum, using the Pieri rule and the following lemma, which paraphrases Theorem 5.4.10 in \cite{MR3287258}.

\begin{lemma}
\label{lem:oops}
Fix a partition $\mu$ with at most $d$ parts, and for a partition $\sigma \in S_d$, let $\sigma \mu$ be the composition defined as
\begin{align*}
(\sigma \mu)_i= \mu_{\sigma^{-1}(i)}-\sigma^{-1}(i)+i. 
\end{align*} 
Let $F(t_1,\dotsb, t_d)$ be a multivariate generating function and let $[\mbf {t}^{\alpha}]F$ be the coefficient of $\mbf{t}^\alpha$, for any composition $\alpha$ with at most $d$ parts. We have
\begin{displaymath}
[\mbf{t}^\mu]\Upsilon F= \sum_{\sigma \in S_d} \text{sgn}(\sigma) \left([\mbf{t}^{\sigma\mu}]F\right).
\end{displaymath}
\end{lemma}

\begin{prop}
\label{prop:altsum}
For a partition $\mu\vdash n$ with at most $d$ parts,
\begin{align}
\label{eq:altsum}
r_{\mu(k|\ell)}= \sum_{j=0}^\ell (-1)^{j}\sum_{\sigma \in S_d} \sgn(\sigma) [u^{k+1+j}v^{\ell-j}\mbf{t}^{\sigma \mu}] \prod_{r\geq 0,1\leq i\leq d}\frac{1+t_iu^rv}{1-t_iu^r}.
\end{align}
\end{prop}

\section{Combinatorial interpretation}
\label{sec:invol}
In this section we interpret Proposition \ref{prop:altsum} combinatorially. We begin by interpreting each summand in Equation \eqref{eq:altsum} as the cardinality of a set of tableau. Let $\ST(\alpha)$ be the set of \emph{standard} tableaus in the alphabet $\mathcal{A}= \{\clb{0}<\clr{0}<\clb{1}<\clr{1}<\dotsb\}$, which satisfy the following properties:
\begin{enumerate}
\item The shape of the tableau is $\alpha$.
\item The entries increase weakly in each row and the red entries increase \emph{strictly} in each row.
\end{enumerate}

For a tableau $T \in \ST(\alpha)$, the \emph{weight} of $T$ is defined as $\wt(T)=\text{ sum of entries of $T$ ignoring their colour}$. Let $\ST(\alpha,\clr{r},{w})$ be the set of standard tableaus of shape $\alpha$ and weight $w$ with $r$ red entries.
\begin{prop}
\label{prop:Hkel}
For a partition $\mu\vdash n$ with at most $d$ parts and a permutation $\sigma \in S_d$ and $0\leq j \leq \ell$, $$[u^{k+1+j}v^{\ell-j}\mbf{t}^{\sigma \mu}]\prod_{r,1\leq i\leq d}\frac{1+t_iu^rv}{1-t_iu^r}=|ST(\sigma \mu, \clr{\ell-j},{k+1+j})|.$$
\end{prop}

\begin{proof}
We will define a bijection between monomials contributing to the coefficient of $u^{k+1+j}v^{\ell-j}\mbf{t}^{\sigma \mu}$ and tableau in $ST(\sigma \mu, \clr{\ell-j},{k+1+j})$. 

Let $T$ be a tableau in $ST(\sigma \mu, \clr{\ell-j},{k+1+j})$. Associate a monomial to a cell $(i,j) \in T$ in the following way: $t_{\sigma^{-1}(i)}u^rv$ if it contains $\clr{r}$ and $t_{\sigma^{-1}(i)}u^r$ if it contains $\clb{r}$ for $r\geq 0$. The monomial associated to $T$ is the product of monomials associated to each cell of $T$. The exponent of $u$ counts the weight of the tableau, which is $k+1+j$ and the exponent of $v$ is the number of red entries, viz. $\ell-j$. The exponent of $t_{\sigma^{-1}(i)}$ is $(\sigma \mu)_{i}$. 

Conversely, we associate a tableau to each monomial in the following way: for each $i$, collect all terms $t_{\sigma^{-1}(i)}u^r$ and $t_{\sigma^{-1}(i)}u^rv$ for any $r\geq 0$. Each of the former type are recorded as $\clb{r}$, while the latter type are recorded as $\clr{r}$. Arrange these numbers in ascending order in the $i$th row of the tableau. There are clearly $(\sigma \mu)_i$ entries in the $i$th row, and the terms are weakly increasing within each row, and the red entries do not repeat within a row. The resulting tableau is therefore in $ST(\sigma \mu, \clr{\ell-j},{k+1+j})$.
\end{proof}

Proposition \ref{prop:Hkel} allows us to write
\begin{align}
\label{eq:hrc_alt}
r_{\mu (k|\ell)}= \sum_{j=0}^l (-1)^j \sum_{\sigma \in S_d}\sgn(\sigma) |ST(\sigma \mu, \clr{\ell-j},{k+1+j})|.
\end{align}
A tableau $T \in ST(\sigma \mu, \clr{\ell-j},{k+1+j})$ is assigned a permutation $\perm(T)=\sigma$, and an \emph{inner sign} $\inn(T)=\sgn(\sigma)$ and an \emph{outer sign} $\outt(T)=(-1)^j$. An involution on a set of signed combinatorial objects is called a \defn{sign-reversing involution} if it either takes objects to themselves (\emph{fixes them}) or to an object with the opposite sign (\emph{cancels them}). In Subsection \ref{subsec:inn} we define an \emph{inner involution} on each $\brac[\Xi]{j}:= \cup_{\sigma\in S_d}ST(\sigma \mu, \clr{\ell-j},{k+1+j})$ that eliminates all tableau with negative inner sign. The fixed points in $\bigcup_{j}\brac[\Xi]{j}$ are acted on by an \emph{outer involution}, defined in Subsection \ref{subsec:out}.

\subsection{Inner involution}
\label{subsec:inn}
Fix a tableau $T \in \brac[\Xi]{j}$. For $i\geq 2$, find the smallest column $j_i$ such that one of the following conditions is true:
\begin{enumerate}
\item $(i,j_i)\in T$ but $(i-1,j_i)\notin T$ or
\item The entry in $(i,j_i)$ is red and strictly less than the entry in $(i-1,j_i)$ or
\item The entry in $(i,j_i)$ is \clb{blue} and less than or equal to the entry in $(i-1,j_i)$.
\end{enumerate}

Let $S=\{j_i|i\geq 2\}$. If $S$ is empty then the algorithm leaves $T$ unchanged. Such a tableau is called a \textbf{fixed tableau}.

If $S$ is not empty, define $j =\min S$ and let $m$ be the largest value such that $j_{m}=j$. We call $(m,j_{m})$ the \textbf{{spot}}. A \textbf{flip} is performed at the spot by the following procedure:
\begin{enumerate}
\item Every row except the $m$th and $m-1$th is unchanged. The cells $(m,j)$ for $j\leq j_m$ are unchanged. The cells $(m-1,j)$ for $j < j_m$ are unchanged. 
\item For each $j\geq j_m$ the cell $(m-1,j)$ is moved (if it exists) to position $(m,j+1)$, while for each $j > j_m$, the cell $(m,j)$ is moved to $(m-1,j-1)$th position.
\end{enumerate}
Note that if a spot exists the flip changes the shape of the tableau. This is clear when either of the last two conditions above are met. When $(i,j_i)\in T$ but $(i-1,j_i)\notin T$ then $(\sigma \mu)_i - (\sigma \mu)_{i-1}\geq 2$, so at least one cell is moved from the $i$th to the $i-1$th row.

\begin{example}[Two rows]
We will demonstrate the algorithm on the tableau in Figure \ref{fig:fig3}. The \emph{spot} is chosen as the first index in the second row satisfying one of the conditions listed above. This location is $(2,2)$. A flip is performed at this spot as in Figure \ref{fig:fig3}
\begin{figure}[h]
\centering
\begin{ytableau}
\clr 0 & 	\clr 1&\clb 1&\clb 2&\clr 2\\
\clr 0 & *(green)\clb 1 & \clb 3
\end{ytableau}
$\rightarrow$
\begin{ytableau}
\clr 0 & \clb 3\\
\clr 0 & *(green)\clb 1 & \clr 1&\clb 1&\clb 2&\clr 2
\end{ytableau}
\caption{The \emph{spot} is coloured {\color{green}green}. Flipping the tableau at the \emph{spot} results in the tableau on the right, which is in $\ST(\sigma(5,3),\clr{4},{10})$, where $\sigma=(12)$.}
\label{fig:fig3}
\end{figure} \end{example}

\begin{example} [Multiple rows]
We will demonstrate the algorithm on the tableau in Figure \ref{fig:fig5}. The following locations satisfy one of the conditions listed above:
\begin{enumerate}
\item $(2,2)$ since the entry in it is \clb{blue} and \emph{less than} the entry in $(1,2)$ .
\item $(3,2)$ since the entry in it is \clb{blue} and \emph{equal} to the entry in $(2,2)$
\item $(4,3)$ since the entry in it is \clb{blue} and \emph{less than} the entry in $(3,3)$.
\end{enumerate}
We pick the leftmost and lowest of these as the \emph{spot}. A flip is performed at the \emph{spot} as in Figure \ref{fig:fig5}. 
\begin{figure}[h]
\centering
\begin{ytableau}
\clb 0 &\clr 1& \clb 1&\clb 2\\
\clr 0& \clb 1& \clb 1&\clb 2\\
\clr 0& *(green)\clb 1& \clr 2\\
\clr 0 & \clb 1 & \clb 2
\end{ytableau}
$\longrightarrow$
\begin{ytableau}
\clb 0 &\clr 1& \clb 1&\clb 2\\
\clr 0& \clr 2\\
\clr 0& *(green)\clb 1&\clb 1& \clb 1&\clb 2\\
\clr 0 & \clb 1 & \clb 2
\end{ytableau}
\caption{Flipping the tableau at the \emph{spot} (coloured {\color{green}green}) results in the tableau on the right, which is in $\ST(\sigma(4,4,3,3),\clr{5},18)$ where $\sigma=(23)$.}
\label{fig:fig5}
\end{figure}
\end{example}

A standard tableau is called a \defn{supertableau} if it satisfies the following properties
\begin{enumerate}
\item Its shape is a partition
\item The \clb{blue} entries increase strictly in each column.
\end{enumerate} 
Let $\SpT(\mu, \clr{\ell-j},{k+1+j})$ denote the set of supertableau of shape $\mu$ and weight $k+1+j$ with ${\ell-j}$ red entries. 

\begin{theorem}
\label{th:isinv_inn}
The procedure described above is a sign-reversing involution on each $\brac[\Xi]{j}$. The set of fixed tableaus in each $\brac[\Xi]{j}$ is $\SpT(\mu, \clr{\ell-j},{k+1+j})$.  
\end{theorem}
\begin{proof}
We first prove the algorithm is sign-reversing. Consider a tableau in $T \in \ST(\sigma \mu, \clr{\ell-j},{a+1+j})$ that is not fixed by the involution. Let the flip occur at $(m,j_m)$ and let $\overline{T}$ be the tableau obtained as a result. The algorithm changes only the shape of the tableau, leaving the number of red entries and the weight unchanged. The algorithm also leaves the length of every row except the $m$th and $m-1$th unchanged.

We claim that $\overline{T}\in \ST(\sigma s_{m-1} \mu, \clr{\ell-j},{a+1+j})$ where $s_{m-1}=(m-1, m)$. Since $\sigma$ and $\sigma s_{m-1}$ agree on all values except $m-1$ and $m$, the size of these rows must remain the same, as is true for $\overline{T}$. A tableau in $ \ST(\sigma s_{m-1} \mu, \clr{\ell-j},{a+1+j})$ must have $(\sigma \mu)_{m-1}+1$ entries in the $m$th row and $(\sigma \mu)_m-1$ entries in the $m-1$th row. There are $\mu_{\sigma^{-1}(m)}-\sigma^{-1}(m)+m$ entries in the $m$th row of $T$ and $\mu_{\sigma^{-1}(m-1)}-\sigma^{-1}(m-1)+(m-1)$ entries in the $(m-1)$th row of $T$. After flipping, there are $(\sigma \mu)_{m-1}-(j_m-1)+j_m$ entries in the $m$th row and $(\sigma \mu)_m -j_m +(j_m-1)$ entries in the $(m-1)$th row of $\overline{T}$. Moreover, it is easy to verify that $\overline{T}$ remains standard.

Next we prove the process is an involution. This amounts to proving that the location of the \emph{spot} is unchanged in $\overline{T}$ as defined above. Since the entry in the $(m-1,j_m)$ position in $\overline{T}$ is $T(m,j_m+1)$ and $T(m,j_m)\leq T(m,j_m+1)$, with equality achieved only if $T(m,j_m)$ is \clb{blue}, the location $(m,j_m)$ is a candidate for the \emph{spot} in $\overline{T}$. The localness of the process ensures that a new candidate for the \emph{spot} can only be created at $(m-1,j_m)$ or $(m+1,s)$ for $s>j_m$. Neither possibility changes the \emph{spot}. 

Finally we characterise the fixed points of this involution, which are precisely the ones without a \emph{spot}. If $\sigma\neq id$, then at least one part of $\sigma \mu$ is larger than the one before it, which creates a \emph{spot}. So the fixed tableau are a subset of $\ST(\mu,\clr{\ell-j},{a+1+j})$. If the columns of a tableau are not weakly increasing in both entries and strongly increasing in the \clb{blue} entries then a \emph{spot} exists. Therefore the fixed tableaus are precisely $\SpT(\mu, \clr{\ell-j},{k+1+j})$.
\end{proof}
\subsection{Outer involution}
\label{subsec:out}
We know from Theorem \ref{th:isinv_inn} that 
\begin{align}
\label{eq:finalsum}
r_{\mu(k|\ell)}=\sum_{j=0}^\ell (-1)^{j}|\SpT(\mu, \clr{\ell-j},{k+1+j})|.
\end{align} 

We now define a sign-reversing involution on $\SpT(\mu):=\cup_{0\leq j \leq \ell}\SpT(\mu,\clr{\ell-j},{k+1+j})$. Let the magnitude of the largest entry in the lowest row be $r$. 
\begin{enumerate}
\item If the largest entry is $\clr{r}$, change $\clr{r}$ to $\clb{r+1}$.

\item If the largest entry is $\clb{r}$, consider the first instance of $\clb{r}$ in this row:
\begin{enumerate}
\item If the entry to its left is less than $\clr{r-1}$:
\begin{enumerate}
\item The tableau is in $\SpT(\mu,\clr{\ell},{k+1})$ then it is \emph{fixed},
\item otherwise change the $\clb{r}$ to $\clr{r-1}$.
\end{enumerate} 
\item If the entry to its left is $\clr{r-1}$ then change the $\clr{r-1}$ to $\clb{r}$.
\end{enumerate}
\end{enumerate}
\begin{example} [Multiple rows]
We will demonstrate the algorithm on the tableaus in Figure \ref{fig:fig6} and \ref{fig:fig7}. The tableau in Figure \ref{fig:fig6} has $\clb{3}$ as its last entry in the last row, and the entry to its left is $\clr{1}$. We may follow Step $2$a here. 
\begin{figure}[h]
\centering
\begin{ytableau}
\clb 0 &\clb 0& \clb 0&\clb 0\\
\clr 0& \clb 1& \clb 1&\clb 1\\
\clr 0& \clr 1& \clb 2\\
\clr 0 & \clr 1 & \clb 3
\end{ytableau}
$\longrightarrow$
\begin{ytableau}
\clb 0 &\clb 0& \clb 0&\clb 0\\
\clr 0& \clb 1& \clb 1&\clb 1\\
\clr 0& \clr 1& \clb 2\\
\clr 0 & \clr 1 & \clr 2
\end{ytableau}
\caption{Step $1 \leftrightarrow$ Step $2$a.}
\label{fig:fig6}
\end{figure}

In Figure \ref{fig:fig7}, it is not possible to change $\clb{3}$ to $\clr{2}$, since this entry exists to its left. We instead change $\clr{2}$ to $\clb{3}$.
\begin{figure}[h]
\centering
\begin{ytableau}
\clb 0 &\clb 0& \clb 0&\clb 0\\
\clr 0& \clb 1& \clb 1&\clb 1\\
\clr 0& \clr 1& \clb 2\\
\clr 0 & \clr 2 & \clb 3
\end{ytableau}
$\longrightarrow$
\begin{ytableau}
\clb 0 &\clb 0& \clb 0&\clb 0\\
\clr 0& \clb 1& \clb 1&\clb 1\\
\clr 0& \clr 1& \clb 2\\
\clr 0 & \clb 3 & \clb 3
\end{ytableau}
\caption{Step $2$a $\leftrightarrow$ Step $2$b.}
\label{fig:fig7}
\end{figure}

\end{example}
\begin{theorem}
\label{th:isinv_out}
Let $\Xi=\cup_{0\leq j \leq \ell}\SpT(\mu,\clr{\ell-j},{k+1+j})$. The procedure described above is a sign-reversing involution on $\Xi$. 
\end{theorem}
\begin{proof}
Let $T\in \SpT(\mu,\clr{\ell-j},{k+1+j})$ be a tableau that is not fixed by the algorithm. If $\clr{r}$ occurs in the last row, increasing it to $\clb{r+1}$ is a supertableau in $\SpT(\mu,\clr{\ell-j-1},{k+1+j+1})$. The largest element in the last row of this tableau is $\clb{r+1}$, which is reduced by the algorithm to $\clr{r}$, yielding $T$.

If the last entry is $\clb{r}$, the entry to the left of the first $\clb{r}$ in this row is less than or equal to $\clr{r-1}$. If the entry is equal to $\clr{r-1}$, increasing the $\clr{r-1}$ to $\clb{r}$ creates a supertableau in $\SpT(\mu,\clr{\ell-j-1},{k+1+j+1})$. If the entry is less than $\clr{r-1}$ then changing the leftmost $\clb{r}$ to $\clr{r-1}$ creates a supertableau in $\SpT(\mu,\clr{\ell-j+1},{k+1+j-1})$. It is easy to see that the process is an involution in both cases.
\end{proof}

The fixed tableau are bijective with the subset of $\SpT(\mu,\clr{\ell+1},{k})$ described below. 

\begin{theorem}[\textbf{hook restriction coefficients}]
\label{th:hrc}
For a partition $\mu\vdash n$ and $k,\ell \geq 0$, the hook restriction coefficient $r_{\mu (k|\ell)}$ is the cardinality of the subset of $\SpT(\mu,\clr{\ell+1},{k})$ whose largest entry in the last row is {red}. 
\end{theorem}

\begin{corollary}
\label{cor:lamlam}
For any partition $\mu$ and integer $n\geq |\mu|+\mu_1$, we have
\begin{align*}
r_{\mu[n](k|\ell)}=\begin{cases}
0 & |\mu|>k+\ell+1 \text{ or } \mu \neq (k|\ell)\\
1 & \mu =(k|\ell).
\end{cases}
\end{align*}
\end{corollary}
\begin{proof}
We maximise the size $\mu$ by first filling the first row of $\mu[n]$ with $\clb{0}$ and then placing $\ell$ {red} entries and the remaining \clb{blue} entries to minimize their contribution to the weight. The {red} entries contribute minimally when in a single column, where they can all be $\clr{0}$; the \clb{blue} entries contribute minimally in a single row, where they can all be $\clr{1}$. Therefore, $\mu=(k|\ell)$ is the largest size of a tableau contributing to $r_{\mu[n](k|\ell)}$. No larger size is possible, and neither is another $\mu$ of size $k+\ell+1$. The only contributing tableau of shape $\mu[n]$ contains $\clb{0}$ in the first row, $\clr{0}$ in the remaining cells of the first column, and $\clb{1}$ in the remaining entries in the second row. 
\end{proof}
\bibliography{bibfile}
\bibliographystyle{plain}
\end{document}

%% file: preamble.tex
\usepackage{euscript}
\usepackage{amsmath,amsthm,amssymb}

\theoremstyle{plain}
\newtheorem{theorem}{Theorem}[section]
\newtheorem{lemma}[theorem]{Lemma}
\newtheorem{corollary}[theorem]{Corollary}

\theoremstyle{definition}
\newtheorem{remark}[theorem]{Remark}
\newtheorem{definition}[theorem]{Definition}
\newtheorem{example}[theorem]{Example}
\newtheorem{exercise}[theorem]{Exercise}

\newtheorem{prop}[theorem]{Proposition}

\theoremstyle{remark}

\renewcommand{\emptyset}{\varnothing}
\newcommand{\ncom}{\newcommand}

\newcommand{\C}{{\mathbb C}}

\newcommand{\GL}{{\mbox{GL}}}

\newcommand{\res}{{\mbox{res}}}

\ncom{\n}[1]{{{ \|{#1}\|     }}}
\ncom{\ns}{\normalsize}
\ncom{\la}{\lambda}
\ncom{\bm}{\boldmath}
\ncom{\noi}{\noindent}
\ncom{\bq}{\begin{equation}}
\ncom{\eq}{\end{equation}}  
\ncom{\beqn}{\begin{eqnarray*}}
\ncom{\eeqn}{\end{eqnarray*}}  
\ncom{\ba}{\begin{array}}
\ncom{\ul}{\underline}   

\ncom{\ea}{\end{array}}
\ncom{\beq}{\begin{eqnarray}}
\ncom{\eeq}{\end{eqnarray}}  
\ncom{\nno}{\nonumber}
\ncom{\hs}{\mbox{\hspace{.25cm}}}
\ncom{\rar}{\rightarrow}
\ncom{\lrar}{\longrightarrow}
\ncom{\Rar}{\Rightarrow}
\ncom{\noin}{\noindent} 
\ncom{\bc}{\begin{center}}
\ncom{\ec}{\end{center}}  
\ncom{\sz}{\scriptsize}   
\ncom{\fpd}{\Phi(\pi^{'})}
\ncom{\fp}{\Phi(\pi) }
\ncom{\nk}{\left< \begin{array}{c}
                       n\\k \end{array} \right>}
\ncom{\nd}{1^{'},2^{'},\cdots,n^{'}}
\ncom{\de}{\bigtriangleup (F_{2n},\leq)}
\ncom{\del}{\bigtriangleup}
\ncom{\cov}{<\!\!\!\!\cdot }
\ncom{\bt}{\begin{theorem}}
\ncom{\bcon}{\begin{con}}
\ncom{\et}{\end{theorem}}
\ncom{\econ}{\end{con}}
\ncom{\bl}{\begin{lemma}}
\ncom{\el}{\end{lemma}}  
\ncom{\bco}{\begin{corollary}} 
\ncom{\ds}{\displaystyle}
\ncom{\eco}{\end{corollary}}   
\ncom{\bp}{\begin{pro}}  
\ncom{\ep}{\end{pro}}    
\ncom{\bex}{\begin{example}}
\ncom{\eex}{\end{example}}  
\ncom{\bexr}{\begin{exercise}}
\ncom{\eexr}{\end{exercise}}  
\ncom{\bprob}{\begin{problems}}
\ncom{\eprob}{\end{problems}}  
\ncom{\bd}{\begin{definition}}
\ncom{\ed}{\end{definition}}  
\ncom{\brm}{\begin{remark}}   
\ncom{\erm}{\end{remark}}     
\ncom{\bal}{\begin{Algorithm}}
\ncom{\eal}{\end{Algorithm}}  
\ncom{\ol}{\overline}
\ncom{\wh}{\widehat} 

\ncom{\pf}{\noi {\bf Proof.  }}
\ncom{\eprf}{\noi {$\Box$}}
\ncom{\be}{\begin{enumerate}} 
\ncom{\ee}{\end{enumerate}}   
\ncom{\seq}{\subseteq}

\ncom{\zr}{\bf\textcolor{red}}
\ncom{\zb}{\bf\textcolor{blue}}
\ncom{\zg}{\bf\textcolor{green}}
\ncom{\zm}{\bf\textcolor{magenta}}

\definecolor{gold}{rgb}{0.85,.66,0}
\definecolor{gb}{rgb}{0, .5,.5}
\definecolor{rb}{rgb}{0.5, 0,.5}
\definecolor{Pink}{rgb}{1,0.75, 0.8}
\ncom{\zgb}{\bf\textcolor{gb}}
\ncom{\zrb}{\bf\textcolor{rb}}

\makeatletter
\newcommand{\single}{\let\CS=\@currsize\renewcommand{\baselinestretch}{1.5}\tiny\CS}
\newcommand{\oneandahalfspacing}{\let\CS=\@currsize\renewcommand{\baselinestretch}{1.5}\tiny\CS}
\newcommand{\doublespacing}{\let\CS=\@currsize\renewcommand{\baselinestretch}{1.6}\tiny\CS}
\newcommand{\double}{\let\CS=\@currsize\renewcommand{\baselinestretch}{3}\tiny\CS}

\def\notarro{{{\hbox{{\hspace*{-.02in}}$\rightarrow$ } }
{\hbox{$\!\!\!\!\!\!\!$}}{\raise 0.2ex 
\hbox{$\scriptscriptstyle{/}$}}}\hspace{.06in}}
\def \*{^{\mbox{$*$}}}

\newcommand\bnota{\begin{nota} }                         
\newcommand\enota{\end{nota} }                           

\newcommand{\sgn}{{\rm sgn\,}}
\newcommand{\beano}{\begin{eqnarray*}}
\newcommand{\eeano}{\end{eqnarray*}}

\newcommand{\defn}[1]{\textbf{#1}} 
\newcommand{\clr}[1]{{\color{red}#1}} 
\newcommand{\clb}[1]{{\color{blue}#1}} 

\newcommand{\Par}{\mathcal{P}}
\newcommand{\GLn}{\GL_n(\C)} 
\newcommand{\WW}{W_\lambda(\C^n)} 
\newcommand{\VV}[1][\mu]{V_{#1}} 
\newcommand{\wt}{\text{wt}} 
\newcommand{\ST}{\text{ST}} 
\newcommand{\SpT}{\text{SpT}} 
\newcommand{\Sym}{\text{Sym}} 
\newcommand{\ch}{\mbox{\Large $ \chi$ }} 
\newcommand{\brac}[2][\beta]{#1^{(#2)}} 
\newcommand{\mch}[2]{\ensuremath{\left(\kern-.3em\left(\genfrac{}{}{0pt}{}{#1}{#2}\right)\kern-.3em\right)}} 
\newcommand{\mbf}[1]{\mathbf{#1}} 
\newcommand{\perm}{\text{perm}} 
\newcommand{\inn}{\text{in}} 
\newcommand{\outt}{\text{out}} 